\documentclass[12pt]{amsart}
\usepackage[utf8]{inputenc}
\usepackage{amssymb,amsmath,amsthm,mathrsfs,url}
\usepackage{multirow, comment}
\usepackage{graphics}

\theoremstyle{remark}

\numberwithin{theorem}{section} \numberwithin{equation}{section}

\newcommand{\FF}{\mathbb{F}_p}

\newcommand{\R}{\mathbb{R}}

\newcommand{\Q}{\mathbb{Q}}

\newcommand{\Z}{\mathbb{Z}}

\DeclareFontFamily{U}{wncy}{}
\DeclareFontShape{U}{wncy}{m}{n}{<->wncyr10}{}
\DeclareSymbolFont{mcy}{U}{wncy}{m}{n}
\DeclareMathSymbol{\Sh}{\mathord}{mcy}{"58}

\newtheorem*{theorem*}{Theorem}

\title{Ranks of elliptic curves and deep neural networks}
\author{Matija Kazalicki}
\address{Department of Mathematics\\ 
	University of Zagreb\\
	Bijeni\v{c}ka cesta 30\\
	10000 Zagreb\\
	Croatia}
\email{matija.kazalicki@math.hr}

\author{Domagoj Vlah}
\address{Department of Applied Mathematics\\
    Faculty of Electrical Engineering and Computing\\
	University of Zagreb\\
	Unska 3\\
	10000 Zagreb\\
	Croatia}
\email{domagoj.vlah@fer.hr}

\begin{document}

\maketitle

\section{Abstract}

Determining the rank of an elliptic curve $E/\Q$ is a hard problem, and in some applications (e.g. when searching for curves of high rank) one has to rely on heuristics aimed at estimating the analytic rank (which is equal to the rank under the Birch and Swinnerton-Dyer conjecture). 

In this paper, we develop rank classification heuristics modeled by deep convolutional neural networks (CNN). Similarly to widely used Mestre-Nagao sums, it takes as an input the conductor of $E$ and a sequence of normalized $a_p$-s (where $a_p=p+1-\#E(\FF)$ if $p$ is a prime of good reduction) in some range ($p<10^k$ for $k=3,4,5$), and tries to predict rank (or detect curves of ``high'' rank). The model has been trained and tested on two datasets: the LMFDB and a custom dataset consisting of elliptic curves with trivial torsion, conductor up to $10^{30}$, and rank up to $10$. For comparison, eight simple neural network models of Mestre-Nagao sums have also been developed.

Experiments showed that CNN performed better than Mestre-Nagao sums on the LMFDB dataset (interestingly neural network that took as an input all Mestre-Nagao sums performed much better than each sum individually), while they were roughly equal on custom made dataset.

\section{Introduction}

Let $E$ be an elliptic curve over $\Q$ with discriminant $\Delta$ and conductor $N$. A celebrated Mordell's Theorem states that the group of rational points $E(\Q)$ is a finitely generated abelian group isomorphic to $E(\Q)_{tors}\times \Z^r$, where $E(\Q)_{tors}$ is the torsion subgroup of $E(\Q)$ and $r$ is the rank of $E(\Q)$. While the possible torsion subgroups were completely classified by Mazur \cite{Mazur77} and easy to compute, the rank of the elliptic curve is a much more mysterious quantity. Not only it is not known what values can rank attain (although this question was already asked by Poincar\'{e} \cite{Poincare}), but also there is no agreement on whether the rank is unbounded or not. Up until recently, the folklore conjecture was that the rank is unbounded which was challenged by the series of papers \cite{Watkins14,Watkins15,Park_et_all} that predict (based on different heuristic models) that there are only finitely many elliptic curves with rank greater than $21$. The current rank record is $28$ by Elkies \cite{Elkies_28}. For more information about the rank records for curves with fixed torsion subgroup as well for rank in families see \cite{Dujella_webpage}.

One of the reasons why it is hard to find curves of high rank is that computationally it is costly to compute the rank of an elliptic curve, mainly due to the fact that we don't know how to efficiently find rational points on elliptic curves. Descent algorithms commonly in use, eventually boil down to the naive point search on some auxiliary curves. One way around this is to use rank heuristics inspired by the Birch and Swinnerton-Dyer conjecture to filter out probable candidates for high rank elliptic curves.

For each prime of good reduction $p\nmid \Delta$, we define $a_p=p+1-\#E(\FF)$. For $p|N$, we set $a_p=0,-1,$ or $1$  if, respectively,  $E$ has additive, split multiplicative or non-split multiplicative reduction at $p$. The $L$-function attached to $E/\Q$ is then defined as an Euler product
$$L_E(s)=\prod_{p|\Delta}\left(1-\frac{a_p}{p^s} \right)^{-1}\prod_{p \nmid \Delta}\left(1-\frac{a_p}{p^s}+\frac{p}{p^{2s}} \right)^{-1},$$
which converges absolutely for $\Re(s) > 3/2$ and extends to an entire function by the Modularity theorem \cite{Wiles,Breuil_Conrad_Diamond_Taylor}. Birch and Swinnerton-Dyer (BSD) conjecture states that the order of vanishing of $L_E(s)$ at $s=1$ (the quantity known as analytic rank) is equal to the rank of $E(\Q)$.

Mestre \cite{Mestre82} and Nagao \cite{Nagao92}, and later others  \cite{Elkies_new_rank_records,Bober}, motivated by BSD conjecture, considered certain sums (see Section \ref{sec:Nagaove sume} for the list of sums examined in this paper) which heuristically should be able to detect curves of high analytic rank. We call all such sums Mestre-Nagao sums. For example, one of these sums (see Section 2 in \cite{Elkies_new_rank_records})
$$\tilde{S_5}(B) =\sum_{\substack{p<B, \\ \textrm{ good  reduction}}} \log\left( \frac{p+1-a_p}{p}\right)$$
has a property that $exp(-\tilde{S_5}(B))$ is the partial product of $L_E(s)$
\begin{equation}\label{eq:partial}
\prod_{\substack{p<B, \\ \textrm{ good  reduction}}}\left( 1-a_p p^{-s}+p^{1-2s}\right)^{-1},
\end{equation}
evaluated at $s=1$ (ignoring the primes of bad reduction). One expects that $\tilde{S_5}(B)$ should be large if $E$ has a large rank since then the partial product should rapidly approach zero. This sum was used in \cite{Elkies_new_rank_records} as a first step in finding rank-record breaking curves with fixed cyclic torsion $\Z/n\Z$ for $n=2,3,\ldots 7$. 

Recently, some new fundamental results have been discovered with the assistance of deep neural networks in topology and representation theory \cite{Davis_Velickovic}, combinatorics \cite{Wagner}, as well as some applications to problems in statistics \cite{Ivek_Vlah21, Ivek_Vlah22}. In number theory, the utility of machine learning methods was shown in \cite{Hee_Lee_Oliver_Arithmetic, Hee_Lee_Oliver_Arithmetic_Pozdnyakov} where the authors, among other things, successfully used logistic regression for classifying elliptic curves of rank zero and one.

In this paper, we investigate a deep learning algorithm for rank classification based on convolutional neural networks (CNN). These networks take as an input the conductor of the elliptic curve together with the sequence of normalized $a_p$-s (i.e. $a_p/\sqrt{p}$) for $p$ in a fixed range and output the rank of the elliptic curve. We compare its performance to that of Mestre-Nagao sums (defined in Section \ref{sec:Nagaove sume}). A priori, it is not clear how to decide on the rank of the elliptic curve based on the value of its Mestre-Nagao sum (see related question (1) in Section 7 of \cite{Elkies_new_rank_records}), so we train simple fully connected neural network to do that task for us. Since the answer critically depends on the conductor of the elliptic curve, these networks, besides Mestre-Nagao sum, take the conductor of the elliptic curve as an input. Training these networks revealed the optimal cutoff of specific Mestre-Nagao sum for rank classification (for cutoffs of $S_5$ see Figure \ref{fig:cutoffs}).

Architecture and the training process of our neural networks are described in Section \ref{sec:arhitektura_mreze}. For training, we used two datasets. One is LMFDB database \cite{LMFDB} which contains $3,824,372$ elliptic curves defined over $\Q$, divided in $2,917,287$ isogeny classes, with rank between $0$ and $5$ and conductor less then $300,000,000$. The other one is a custom made dataset which consists of $2,033,965$ elliptic curves with trivial torsion, the conductor less than $10^{30}$ and the rank in the range between $0$ and $10$. For more information about datasets see Section \ref{sec:podaci}.

For each neural network (CNN or one of the Mestre-Nagao sums, in total $9$) we have performed $24$ tests by varying
\begin{itemize}
	\item[a)] dataset - LMFDB or custom,
	\item[b)] range of $a_p$-s - we considered $a_p$-s for $p<10^3, 10^4,$ and  $10^5$, 
	\item[c)] test curves - uniformly selected  ($20\%$ from the dataset) or all curves in top conductor range (which is $[10^8,10^9]$ for LMFDB and $[10^{29},10^{30}]$ for custom dataset),
	\item[d)] type of classification - binary or all ranks (for LMFDB the rank range is from $0$ to $5$, and for custom dataset from $0$ to $10$).
\end{itemize}
In binary classification, curves are labeled as either of low or high rank. For LMFDB high rank means rank $4$ (we did not consider $19$ rank $5$ curves), while for custom dataset high rank is $8,9$ or $10$. Here, the idea is to test how good neural networks are in detecting elliptic curves of ``high'' rank. 

The neural networks have been selected to maximize Matthews correlation coefficient (MCC) or phi coefficient, which is generally regarded as a balanced measure of the quality of the classification even if the classes are of very different sizes \cite{Boughorbel2017}. For binary classification, MCC is computed using the following formula
$$MCC=\frac{TP \cdot TN - FP\cdot FN}{\sqrt{(TP+FP)\cdot(TP+FN)\cdot(TN+EP)\cdot(TN+FN)}},$$
where $TP, FN, TN, FP$ denote, respectively, the number of true positives, false negatives, true negatives, and false positives. Note that $MCC$ lies in the segment $[-1,1]$ and $MCC=1$ only in the case of perfect classification.

For comparison of Mathews correlation coefficients of CNN vs. $S_5$ classifiers for uniform test set see Table \ref{tab:CNNvsS5}.

\begin{table}[]
\centering
\begin{tabular}{c|cccccc|}
\cline{2-7}
                                                                                                    & \multicolumn{6}{c|}{Number of $a_p$-s used}                                                                                                                              \\ \hline
\multicolumn{1}{|c|}{\multirow{2}{*}{\begin{tabular}[c]{@{}c@{}}Type of\\ classifier\end{tabular}}} & \multicolumn{3}{c|}{LMFDB}                                                                    & \multicolumn{3}{c|}{custom dataset}                                      \\ \cline{2-7} 
\multicolumn{1}{|c|}{}                                                                              & \multicolumn{1}{c|}{$p<10^3$} & \multicolumn{1}{c|}{$p<10^4$} & \multicolumn{1}{c|}{$p<10^5$} & \multicolumn{1}{c|}{$p<10^3$} & \multicolumn{1}{c|}{$p<10^4$} & $p<10^5$ \\ \hline
\multicolumn{1}{|c|}{CNN}                                                                           & \multicolumn{1}{c|}{0.9507}   & \multicolumn{1}{c|}{0.9958}   & \multicolumn{1}{c|}{0.9992}   & \multicolumn{1}{c|}{0.6129}   & \multicolumn{1}{c|}{0.7218}   & 0.7958   \\ \hline
\multicolumn{1}{|c|}{$S_5$}                                                                         & \multicolumn{1}{c|}{0.6132}   & \multicolumn{1}{c|}{0.7774}   & \multicolumn{1}{c|}{0.8463}   & \multicolumn{1}{c|}{0.4987}   & \multicolumn{1}{c|}{0.5990}   & 0.6696   \\ \hline
\end{tabular}
\vskip 0.1cm
\caption{Comparison of Mathews correlation coefficients of CNN vs. $S_5$ classifiers for uniform test set.}
\label{tab:CNNvsS5}
\end{table}

\begin{figure}
    \centering
    \resizebox{140mm}{!}{\includegraphics{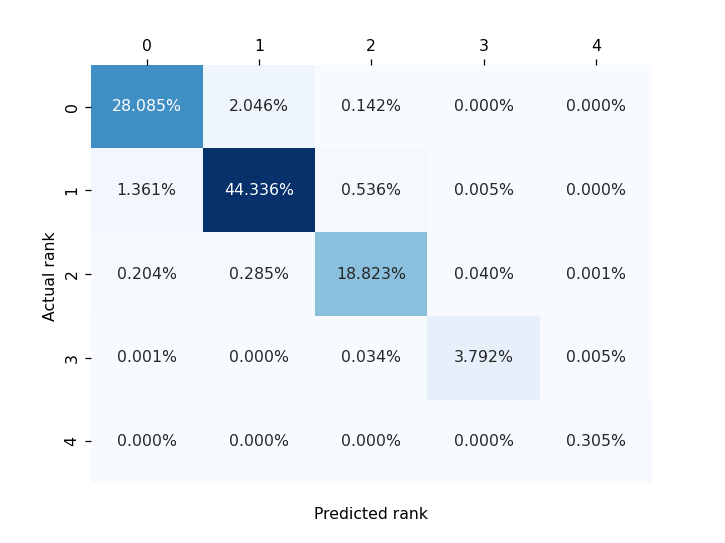}}
    \caption{Confusion matrix of CNN for LMFDB and $p<10^4$ with top conductor range.}
    \label{fig:CM_LMFDB_range_NN_p10000}
\end{figure}

\begin{figure}
    \centering
    \resizebox{140mm}{!}{\includegraphics{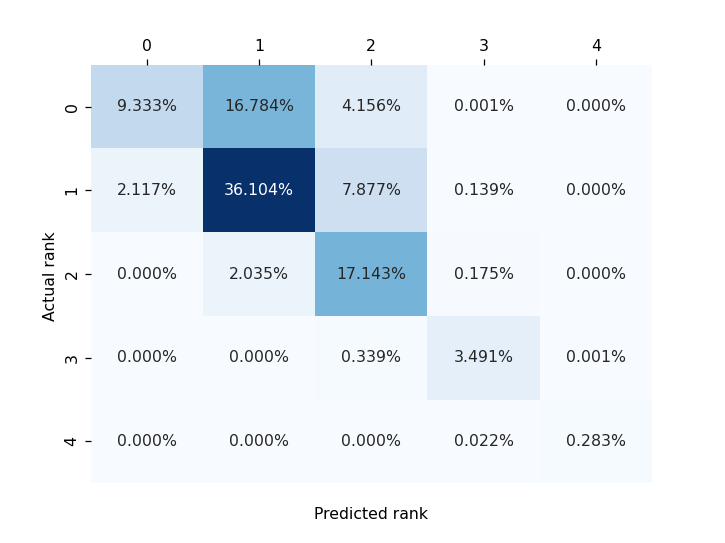}}
    \caption{Confusion matrix of $S_0$ for LMFDB and $p<10^4$ with top conductor range.}
    \label{fig:CM_LMFDB_range_S0_p10000}
\end{figure}

In all the tests performed on LMFDB dataset CNN outperformed Mestre-Nagao sums (this is especially true in the classification of all ranks). For example, in the all ranks classification with $p<10,000$ in uniform range MCC of CNN is $0.9958$ (in particular, CNN misclassified only $0.25\%$ of the curves) while the best Mestre-Nagao sum $S_2$ has $CNN=0.8697$ (it misclassified  $8.1\%$ of the curves). For more details see Table \ref{tab:LMFDB-uniform}. Similarly, in the all ranks classification with $p<10,000$ and top conductor range, MCC of CNN is $0.9289$ (see confusion matrix in Figure \ref{fig:CM_LMFDB_range_NN_p10000}) while the best Mestre-Nagao sum $S_0$ has $MCC=0.5057$ (see Figure \ref{fig:CM_LMFDB_range_S0_p10000}), which is remarkable!

We also trained a fully connected neural network $\Omega$ whose inputs were all Mestre-Nagao sums \{$S_0,\ldots,S_6$\} (together with the conductor), and interestingly on the LMFDB dataset this network performed much better than any sum individually. For example, in the uniform range case above it attained $MCC=0.9602$ and in the top range case $MCC=0.7013$.

A binary classification on the LMFDB dataset (which tries to identify curves of rank $4$) ended up being an easy task for both CNN and the majority of Mestre-Nagao sums - in a uniform mode they demonstrated the prefect classification even in $p<1,000$ range (see Table \ref{tab:LMFDB-uniform} and \ref{tab:LMFDB-range}).

The classification on the custom dataset was much more challenging for both CNN and Mestre-Nagao sums. For example, in all ranks classification with $p<10,000$ and uniform range, while CNN with $MCC=0.7218$ (it misclassified $23\%$ of curves while for $3\%$ of the curves prediction missed true rank for more than $1$) still outperformed the best Mestre-Nagao sum $S_1$ with $MCC=0.6890$ (it misclassified $26\%$ of curves), it performed just a little bit better than the network $\Omega$ for which $MCC=0.7069$ (see Table \ref{tab:triv_5-uniform}).

In the top conductor range (for all ranks classification on the custom dataset) all heuristics were much less efficient. This is expected, at least for Mestre-Nagao sums, since the curves in this test set have large conductors (at least $10^{29}$) and thus we need more $a_p$'s to determine their $L$-function. Mestre-Nagao sums slightly outperformed CNN in this setting. For example, in $p<10,000$ mode MCC of CNN was $0.3019$ (it misclassified $61\%$ of curves while for $12\%$ of the curves prediction missed true rank for more than $1$), while the best Mestre-Nagao sum $S_2$ had $MCC=0.3291$. For complete information see Table $\ref{tab:triv_5-range}$.

In binary classification on the custom dataset (which tries to identify curves of rank $8,9$ or $10$), CNN outperformed each Mestre-Nagao sum individually, but it was a little bit worse than the neural network trained on all Mestre-Nagao sums. For example, in the top conductor range for $p<10,000$ MCC of CNN was $0.6774$ while MCC of all sums was $0.7091$.

\section{Mestre-Nagao sums}\label{sec:Nagaove sume}
In this section, we define and motivate Mestre-Nagao sums that we are going to use for rank classification.

Let $\Lambda(n)=\begin{cases}
	\log{p} &\textrm{if $n=p^m$}\\
	0 &\textrm{otherwise}
\end{cases}$
be the von Mangoldt function, and for elliptic curve $E$ defined over $\Q$ with conductor $N_E$ let
$$
c_n=\begin{cases}
	\alpha_p^m+\beta_p^m, &\textrm{if $n=p^m$ and $p\nmid N_E$},\\
	a_p^m, &\textrm{ if $n=p^m$ and $p|N_E$},\\
	0, &\textrm{otherwise},
\end{cases}
$$
where $\alpha_p$ and $\beta_p$, for $p \nmid N_E$, are eigenvalues of the Frobenius morphism at $p$ (thus $a_p=\alpha_p+\beta_p$ and $\alpha_p\beta_p=p$).
Note that for $m>1$ we have $$\alpha_p^m+\beta_p^m=(\alpha_p^{m-1}+\beta_p^{m-1})a_p-p(\alpha_p^{m-2}+\beta_p^{m-2}),$$
hence for $p\nmid N_E$ it follows $c_{p^m}=c_{p^{m-1}}a_p-p c_{p^{m-1}}$. These number are related to analytic rank $r$ of elliptic curve $E/\Q$ by the following identity 
$$-\frac{L'_E(s)}{L_E(s)}=\sum_{n=1}^\infty \frac{c_n \Lambda(n)}{n^s}=\frac{r}{s-1}+\cdots.$$

We consider the following Mestre - Nagao sums.
\begin{align*}
 S_0(B)&=\frac{1}{\log{B}}\sum_{\substack{p<B, \\ \textrm{ good  reduction}}} \frac{a_p(E)\log{p}}{p},\\
 S_1(B)&=S_0(B)-\frac{1}{B \log{B}}\sum_{n \le B} c_n \Lambda(n),\\
 S_2(B)&=\frac{1}{\log{B}}\sum_{n \le B} \frac{c_n \Lambda(n)}{n}-\frac{1}{B \log{B}}\sum_{n \le B}c_n \Lambda(n),\\
 S_3(B)&= \sum_{\substack{p<B, \\ \textrm{ good  reduction}}} \frac{-a_p(E)+2}{p+1-a_p(E)}\log{p},\\
 S_4(B)& = \frac{1}{B} \sum_{\substack{p<B, \\ \textrm{ good  reduction}}} -a_p(E)\log{p},\\
 S_5(B)&=\sum_{\substack{p<B, \\ \textrm{ good  reduction}}} \log\left( \frac{p+1-a_p(E)}{p}\right)+\sum_{\substack{p<B, \\ \textrm{ split mult. reduction}}} \log\left( 3/2 \cdot \frac{p-1}{p}\right),\\
 S_6(\Delta) &=\frac{\log{N_E}}{2 \Delta \pi}-\frac{\log{2\pi}}{\Delta \pi}+\frac{-1}{\Delta \pi}\sum_{p \le \exp(2\pi \Delta)} \log{p}\sum_{k=1}^{\lfloor 2\pi \Delta/\log{p}\rfloor}\frac{c_{p^k}}{p^{k/2}}\left(1-\frac{k \log{p}}{2\pi \Delta}\right)\\
 &+\frac{1}{\pi}\Re\left\{\int_{-\infty}^{\infty}\frac{\Gamma'}{\Gamma}(1+i t)\left(\frac{\sin(\Delta\pi t)}{\Delta \pi t}\right)^2 dt \right\} .
\end{align*}

The sum $S_0$ was analyzed in detail in \cite{Kim_Murty}, where it was shown that if the Riemann hypothesis for $L_E(s)$ is true, and if there exist the limit $\lim_{B\rightarrow \infty} S_0(B)$, then the limit is $-r+1/2$ where $r$ is analytic rank of $E/\Q$. The relation between $S_0$ and analytic rank comes from the Perron's formula applied to $$\frac{1}{2\pi i}\int -\frac{L_E'(s)}{L_E(s)}\frac{x^s}{s(s-1)}ds, \quad x>1,\, d>3/2,\, a\in \R.$$
From the proof of Theorem 6. in \cite{Kim_Murty}, by modifying the error term of $S_0(B)$ one obtains sums $S_1$ and $S_2$ (equations (3.11) and (3.12)) that have essentially the same limit as $S_0$, but faster convergence.

Sums $S_3$ and $S_4$ were considered by Nagao \cite{Nagao92}. Note that $S_3(B)$ is logarithmic derivative of the partial Euler product \eqref{eq:partial} evaluated at $s=1$ (ignoring the primes of bad reduction). Since $\frac{L_E'(s)}{L_E(s)}=\frac{r}{s-1}+\cdots$, it is reasonable to assume that for elliptic curve with large analytic rank $r$ the corresponding $S_3$ sum will be large (see Section 1.3.3 in \cite{Campbell_disertation}).
Sum $S_5$ is the modification of $\tilde{S_5}$ (from Introduction) at primes of split multiplicative reduction (see Section 7. in \cite{Elkies_new_rank_records}). 

Sum $S_6$ is described in \cite{Bober}. It is based on the explicit formula for $L_E(s)$ relating $S_6(\Delta)$ to the sum over nontrivial zeros $\frac{1}{2}+i\gamma$ of $L_E(s)$, $\sum_\gamma f(\gamma)$, where $f(z)=f(z;\Delta)=\left(\frac{\sin(\Delta\pi z)}{\Delta \pi z}\right)^2$. Note that $f(0)=1$ and $f(x)\ge 0$ for all real $x$, thus, $\sum_\gamma f(\gamma)$ gives an upper bound for $r$. Assuming Riemann Hypothesis for $L_E(s)$, one has $\lim_{\Delta\rightarrow \infty} S_6(\Delta)=r$, but computing $S_6(\Delta)$ becomes infeasible once $\Delta$ gets a little larger than $4.4$. For example, already to compute $S_6(2)$ one needs $a_p$'s for all $p<286 ,751$. 

For computing $S_6(\Delta)$ we used C version \cite{smalljac} of PSAGE \cite{psage}  rankbound algorithm by Bober.

\section{Datasets} \label{sec:podaci}

We have been training our models on two datasets. 
\subsection{LMFDB} LMFDB database \cite{LMFDB} contains $3,824,372$ elliptic curves defined over $\Q$, distributed in $2,917,287$ isogeny classes. We selected one representative per isogeny class (since isogeny curves have equal $L$-functions and consequently same analytic rank). Since there are only $19$ curves with rank $5$, we didn't take any on them in consideration. See Table \ref{tab:rank_distribution_LMFDB} for rank distribution.
Moreover, $63.71\%$ of curves have trivial torsion, while $31.18\%$ have $\Z/2\Z$ torsion.
\begin{table}
	\begin{tabular}{|l|l|l|l|l|l|l|}
		\hline
		rank       & 0       & 1       & 2       & 3      & 4      & 5      \\ \hline
		count      & 1,404,510 & 1,887,132 & 493,291  & 37,334  & 2,086   & 19     \\ \hline
		proportion & 36.73\% & 49.34\% & 12.90\% & 0.98\% & 0.05\% & 0.00\% \\ \hline
	\end{tabular}
	\vskip 0.2cm
	\caption{Distribution of ranks in LMFDB.}
	\label{tab:rank_distribution_LMFDB}
\end{table}

\subsection{Custom dataset}
Our custom dataset contains $2,074,863$ elliptic curves with trivial torsion and conductor less than $10^{30}$. See Table \ref{tab:conductor_distribution} for the number of curves in a specific conductor range. Curves have a rank between $0$ and $10$.
See Table \ref{tab:rank_distribution_custom} for distribution of curves by rank.

\begin{table}[]
\centering
\begin{tabular}{|l|l|l|l|l|l|l|l|}
\hline
range      & $1,10^5$ & $10^5, 10^{10}$ & $10^{10}, 10^{15}$ & $10^{15}, 10^{20}$ & $10^{20}, 10^{25}$ & $10^{25}, 10^{29}$ & $10^{29}, 10^{30}$ \\ \hline
proportion & 0.94\%     & 21.60\%           & 26.97\%              & 21.27\%              & 17.02\%              & 10.19\%              & 1.97\%               \\ \hline
\end{tabular}
\caption{Distribution of conductors in custom dataset.}
\label{tab:conductor_distribution}
\end{table}

\begin{table}
\centering
\begin{tabular}{|l|l|l|l|l|l|l}
\hline
rank       & 0         & 1         & 2         & 3         & 4         & \multicolumn{1}{l|}{5}         \\ \hline
count      & $78,755$  & $121,966$ & $179,593$ & $435,195$ & $543,713$ & \multicolumn{1}{l|}{$405,129$} \\ \hline
proportion & 3.79\%    & 5.87\%    & 8.65\%    & 20.97\%   & 26.20\%   & \multicolumn{1}{l|}{19.52\%}   \\ \hline \hline
rank       & 6         & 7         & 8         & 9         & 10        &                                \\ \cline{1-6}
count      & $210,439$ & $77,923$  & $19,272$  & $2,670$   & $201$     &                                \\ \cline{1-6}
proportion & 10.14\%   & 3.75\%    & 0.92\%    & 0.12\%    & 0.01\%    &                                \\ \cline{1-6}
\end{tabular}
\caption{Distribution of ranks in custom dataset.}
\label{tab:rank_distribution_custom}
\end{table}

We used two methods for generating the dataset. First, we produced curves with random Weierstrass coefficients. These were mostly curves of rank $0$ and $1$. 

The curves of higher rank were obtained as random specializations of pencils of cubics through randomly selected $k$ rational points in the plane, for $k=2,3,\ldots, 8$. Every such cubic (if nonsingular) was transformed to Weierstrass form by the change of variables described in Chapter 8 of \cite{Cassels}. We removed curves with nontrivial torsion subgroup.

Next, using PARI/GP \cite{Pari} function ellrank we tried to compute ranks of all previously generated curves (assuming the Parity conjecture). In the end, we discarded those for which PARI/GP couldn't find the rank. Note that ellrank in computing the upper bound for the rank of the elliptic curve $E$ approximates the rank of $\Sh(E/\Q)[2]$ from below with the rank of $\Sh(E/\Q)[2]/2\Sh(E/\Q)[4]$ (which it computes via Cassels pairing), thus it can not determine the rank of $E$ if $\Sh(E)[4]$ is nontrivial. Consequently, such curves are missing from our dataset.

\section{Neural network architectures and training procedure}
\label{sec:arhitektura_mreze}

The models we employ are convolutional neural networks (CNN). CNN could be regarded as the composition of functions, so-called layers. Each layer is in principle composition of a fittable affine mapping and an element-wise nonlinear function. To detect and utilize any correlations in $a_p$ sequence inherent to our data, for the affine mapping we use 1D matrix convolutional operations, inspired by the success of 2D matrix convolutions in image classification task \cite{Krizhevsky2012}. The point is that all Mestre-Nagao sums inherently employ the same computation ``rule'' for each $a_p$, which we expect could be efficiently modeled by a matrix convolution operation. Thus, we can regard the training of a CNN as an optimization in the highly dimensional parameterized space of Mestre-Nagao like heuristics.

Additionally, we use smaller fully connected neural networks (FCNN) for the classification of curves using only precomputed values of selected Nagao-Mestre sums and conductors.

All described models are implemented in Python using the PyTorch library \cite{Paszke2019}. 
Our models training and evaluation code is available online \cite{githubrepo}.
Models were trained using a server with dual Xeon 5218R CPUs (40 cores) equipped with 6 Nvidia Quadro RTX 5000 GPUs, each with 16\,GB of RAM.

\subsection{Architecture of CNNs}

For each elliptic curve $E$, CNN takes as an input a matrix $M(E)$ of shape $3\times\pi(N)$, where $N=10^3, 10^4, 10^5$ and the columns are indexed with primes less than $N$. First row consist of $a_p$-s divided by $\sqrt{p}$ (by Hasse theorem $a_p/\sqrt{p} \in[-2,2]$). Second row contains constant sequence of normalized conductor, $\log{N_E}/\log{N_{max}}$, where $N_{max}$ is the maximal conductor present in the dataset. Third row is a sequence of $r_n$ (linear sweep from $-1$ do $1$) that for $n$-th prime $p_n$ is equal to $-1 + 2 n/\pi(N)$. In this way we input to the network the relative position of prime $p$ in sequence of primes less than $N$.

Each CNN consist of the following sequence of convolutional layers:
\begin{enumerate}
    \item preparatory layer (increasing number of channels (rows in input matrix) from $3$ to $64$)
    \item several input layers $I_1,\ldots,I_{L_1}$, preserving the number of channels and the length of the input sequence
    \item several reducing layers $R_1,\ldots,R_{L_2}$, preserving the number of channels and reducing the length of input sequence by a factor $2$ (convolution kernel stride is equal to $2$)
    \item several output layers $O_1,\ldots,O_{L_3}$, preserving the number of channels and the length of the input sequence
\end{enumerate}
Each convolution layer is followed by an activation function and a batch normalization layer, improving the convergence of the network training and network generalization properties \cite{Ioffe2015}. For the activation function, we decided to use ReLU. Kernel size $KS$ is the same in each layer. The number of reducing layers $L_2$ is chosen such that the length of the output of the final reducing layer $R_{L_2}$ (and consequently the length of the output of every output layer) is equal to $1$.

Following convolutional layers, we have one final fully connected layer reducing dimensionality from $64$ to the number of classification labels (ranks). Thus, our CNN architecture is fully determined with hyperparameters $L_1$, $L_2$, $L_3$, and $KS$.

For instance, in the case of uniform range for LMFDB model with $p<10,000$ we used the following values of hyperparameters: $L_1=0$, $L_2=11$, $L_3=3$, $KS=17$. To find optimal values of hyperparameters we conducted Bayesian optimization \cite{Mockus1975} in the space of hyperparameters by training altogether $500$ different variants of networks. Similar hyperparameters were acquired for LMFDB in the case of $p<1,000$ and for $p<100,000$, due to lengthy computations we just reused the hyperparameters for $p<10,000$, by setting $L_2=14$. For each custom dataset CNN model, we took the same values of hyperparameters as for the corresponding LMFDB model.

\subsection{Architecture of FCNNs} Here the input to the network for each curve $E$ consists of one or more precomputed values of Nagao-Mestre sums for $a_p$ in a certain range and normalized conductor $N_E$, as in the case of CNNs. In every case, our FCNN consists of an input linear layer having the size of the number of input features (Nagao-Mestre sums and normalized conductor), four hidden linear layers of size $128$, and an output linear layer of the size of the number of classification labels (ranks). In between each of the linear layers, we use dropout \cite{Nitish2014}, which reduces overfitting during training, and ReLU activation function. In our experience, the capacity of this network is more than enough to learn the dependency of rank to the input.

\subsection{Training procedure}
To train a neural network means to optimize network parameters (weights in convolutional and linear layers) in order to minimize a chosen loss function over a given dataset.

\subsubsection*{Dataset}
For neural network training, as customary the whole dataset is divided into three disjoint subsets called training, validation, and test dataset. Only the training and validation datasets are used in network training, thus the test dataset is not ''seen`` by the network during the training and is only used at the end to test the performance of already trained models. We used the $4:1$ ratio between the training and validation datasets.

\subsubsection*{Loss function}
As we are training classifier networks, we decided to use the cross entropy loss function as implemented in PyTorch library \cite{Paszke2019}. Additionally, we supply weights for the loss function to reflect on disproportional sizes of classes, which is considered standard practice.

\subsubsection*{Optimizer hyperparameters}
\label{sec:hyper}

Our models are trained using the fast.ai library \cite{Howard2018}. For the optimizer, we use the usual choice of Adam algorithm \cite{kingma2017adam}. For governing learning rate and momentum scheduling we used a one-cycle policy with cosine annealing \cite{smith2018superconvergence}, which ensures better convergence of model parameters to a broad optimum and allows for better generalization of the trained model \cite{smith2018disciplined}. For the optimizer hyperparameters, we used different ranges of the maximum learning rate for CNNs and FCNNs.

For CNNs the range was between $3.5\cdot 10^{-4}$ and $1.5\cdot 10^{-3}$. Number of training epochs was equal to $40$, and batch size equal to $1024$ ($p<10^5$) or $2048$ ($p<10^3$ and $p<10^4$). For Adam-specific hyperparameters, we have used: $\beta_2=0.99$, $\epsilon=10^{-5}$, and weight decay of $10^{-3}$.

For FCNNs, we employed Bayesian optimization for three hyperparameters: dropout, maximum learning rate, and weight decay, by training $200$ different models (networks) for each of $7$ Mestre-Nagao sums $\{S_0,\ldots,S_6\}$ and also for $\Omega$, in each of the $24$ test variants. Finally, in each of these $192$ tests performed we took the model having the best MCC on the test set. The number of training epochs for every model was equal to $10$ and the rest of the Adam-specific hyperparameters were the same as in CNN training.

\section{Results}

We present rank classification results using the Mathews correlation coefficient in each of the $216$ tests. We use two different datasets (LMFDB and custom, see Section \ref{sec:podaci}), three $a_p$-s ranges ($p<10^3, 10^4, 10^5$), two modes for selecting test curves ($20\%$ uniformly selected or top conductor range), and two types of classification (binary or all ranks).

For binary classification we always trained our networks for all ranks classification and to obtain binary classification we merged all rank labels into two categories, ''low`` and ''high`` ranks - effectively only reinterpreting all ranks classification, but rerunning  Bayesian optimization for hyperparameters. For LMFDB low ranks are $0,1,2,3$ and high ranks are $4$, and for custom dataset low ranks are $0,\ldots,7$ and high ranks are $8,9,10$. Interestingly, this procedure produced much better binary classification results than training directly a binary classifier model.

Regarding CNN models, for each out of $12$ test variants (varying two datasets, three ranges, and two modes for selecting test curves) we trained only $10$ different networks due to high computational expense. For instance training of models in the worst case of $p<10^5$ took over $14$ hours per model on a single GPU. On the other hand, training of FCNNs for Mestre-Nagao sums was much less computational demanding, so for each of the $192$ test cases, we trained at least $200$ models (which took several minutes per model on a single GPU), using Bayesian optimization in the space of hyperparameters, as described in Section \ref{sec:hyper}. For each of the test cases, we selected the best-performing model out of all trained models.

\subsection{LMFDB}

Test results are presented in Tables \ref{tab:LMFDB-uniform} and \ref{tab:LMFDB-range}.

In Table \ref{tab:LMFDB-uniform} are maximal values of MCC obtained for different classifiers. Classifiers are trained using $80\%$ of randomly selected curves from the LMFDB dataset, having a conductor less than $10^8$. MCC is computed on the other $20\%$ of the dataset, not seen during the training. Classifiers (CNN, $\Omega$, and Mestre-Nagao sums from Section \ref{sec:Nagaove sume}) are trained using values of $a_p$-s for primes $p$ less than $10^3$, $10^4$ or $10^5$, and the value of the elliptic curve conductor.

\begin{table}
\centering
\begin{tabular}{c|cccccc|}
\cline{2-7}
                                                                                                    & \multicolumn{6}{c|}{Number of $a_p$-s used}                                                                                                                                                   \\ \hline
\multicolumn{1}{|c|}{\multirow{2}{*}{\begin{tabular}[c]{@{}c@{}}Type of\\ classifier\end{tabular}}} & \multicolumn{3}{c|}{classify all ranks}                                                                            & \multicolumn{3}{c|}{binary classification}                               \\ \cline{2-7} 
\multicolumn{1}{|c|}{}                                                                              & \multicolumn{1}{c|}{$p<10^3$}        & \multicolumn{1}{c|}{$p<10^4$}        & \multicolumn{1}{c|}{$p<10^5$}        & \multicolumn{1}{c|}{$p<10^3$} & \multicolumn{1}{c|}{$p<10^4$} & $p<10^5$ \\ \hline
\multicolumn{1}{|c|}{CNN}                                                                           & \multicolumn{1}{c|}{\textbf{0.9507}} & \multicolumn{1}{c|}{\textbf{0.9958}} & \multicolumn{1}{c|}{\textbf{0.9992}} & \multicolumn{1}{c|}{1}        & \multicolumn{1}{c|}{1}        & 1        \\ \hline
\multicolumn{1}{|c|}{$S_0$}                                                                         & \multicolumn{1}{c|}{0.6823}          & \multicolumn{1}{c|}{0.8435}          & \multicolumn{1}{c|}{0.9068}          & \multicolumn{1}{c|}{0.9877}   & \multicolumn{1}{c|}{0.9758}   & 0.9729   \\ \hline
\multicolumn{1}{|c|}{$S_1$}                                                                         & \multicolumn{1}{c|}{0.6848}          & \multicolumn{1}{c|}{0.8507}          & \multicolumn{1}{c|}{0.9301}          & \multicolumn{1}{c|}{0.9968}   & \multicolumn{1}{c|}{1}        & 1        \\ \hline
\multicolumn{1}{|c|}{$S_2$}                                                                         & \multicolumn{1}{c|}{0.7277}          & \multicolumn{1}{c|}{\textbf{0.8697}} & \multicolumn{1}{c|}{0.9359}          & \multicolumn{1}{c|}{0.9938}   & \multicolumn{1}{c|}{1}        & 1        \\ \hline
\multicolumn{1}{|c|}{$S_3$}                                                                         & \multicolumn{1}{c|}{0.6933}          & \multicolumn{1}{c|}{0.8499}          & \multicolumn{1}{c|}{0.9142}          & \multicolumn{1}{c|}{0.9420}   & \multicolumn{1}{c|}{0.9842}   & 0.9726   \\ \hline
\multicolumn{1}{|c|}{$S_4$}                                                                         & \multicolumn{1}{c|}{0.2678}          & \multicolumn{1}{c|}{0.3015}          & \multicolumn{1}{c|}{0.1525}          & \multicolumn{1}{c|}{0.2188}   & \multicolumn{1}{c|}{0.1103}   & 0.0744   \\ \hline
\multicolumn{1}{|c|}{$S_5$}                                                                         & \multicolumn{1}{c|}{0.6132}          & \multicolumn{1}{c|}{0.7774}          & \multicolumn{1}{c|}{0.8463}          & \multicolumn{1}{c|}{0.9968}   & \multicolumn{1}{c|}{0.9968}   & 1        \\ \hline
\multicolumn{1}{|c|}{$S_6$}                                                                         & \multicolumn{1}{c|}{\textbf{0.6969}} & \multicolumn{1}{c|}{0.8647}          & \multicolumn{1}{c|}{\textbf{0.9381}} & \multicolumn{1}{c|}{1}        & \multicolumn{1}{c|}{1}        & 1        \\ \hline
\multicolumn{1}{|c|}{$\Omega$}                                                                      & \multicolumn{1}{c|}{\textbf{0.8685}} & \multicolumn{1}{c|}{\textbf{0.9602}} & \multicolumn{1}{c|}{\textbf{0.9826}} & \multicolumn{1}{c|}{1}        & \multicolumn{1}{c|}{1}        & 1        \\ \hline
\end{tabular}

\caption{LMFDB with uniform test dataset.}
\label{tab:LMFDB-uniform}
\end{table}

\begin{figure}
    \centering
    \resizebox{140mm}{!}{\includegraphics{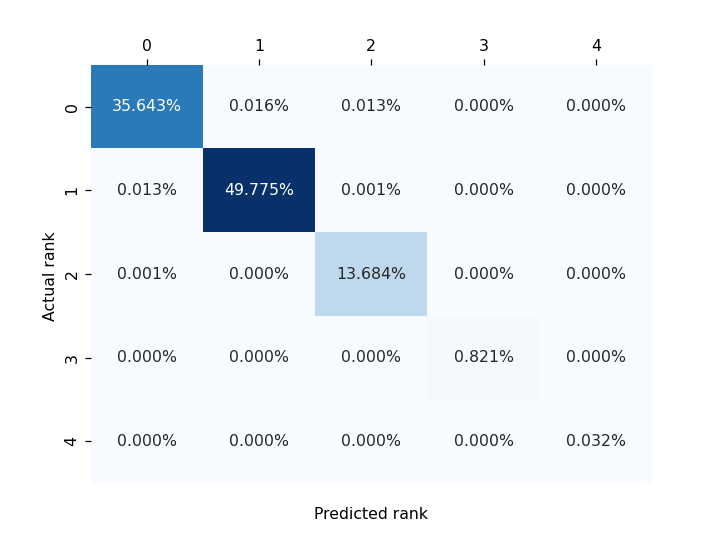}}
    \caption{Confusion matrix of CNN for LMFDB and $p<10^5$ with uniform test dataset.}
    \label{fig:CM_LMFDB_uniform_NN_p100000}
\end{figure}

In Table \ref{tab:LMFDB-range} classifiers are trained using all curves from the LMFDB dataset having a conductor less than $10^8$. MCC is computed on all curves from the LMFDB dataset having a conductor between $10^8$ and $10^9$, not seen during the training.

\begin{table}
\centering
\begin{tabular}{c|cccccc|}
\cline{2-7}
                                                                                                    & \multicolumn{6}{c|}{Number of $a_p$-s used}                                                                                                                                                   \\ \hline
\multicolumn{1}{|c|}{\multirow{2}{*}{\begin{tabular}[c]{@{}c@{}}Type of\\ classifier\end{tabular}}} & \multicolumn{3}{c|}{classify all ranks}                                                                            & \multicolumn{3}{c|}{binary classification}                               \\ \cline{2-7} 
\multicolumn{1}{|c|}{}                                                                              & \multicolumn{1}{c|}{$p<10^3$}        & \multicolumn{1}{c|}{$p<10^4$}        & \multicolumn{1}{c|}{$p<10^5$}        & \multicolumn{1}{c|}{$p<10^3$} & \multicolumn{1}{c|}{$p<10^4$} & $p<10^5$ \\ \hline
\multicolumn{1}{|c|}{CNN}                                                                           & \multicolumn{1}{c|}{\textbf{0.5631}} & \multicolumn{1}{c|}{\textbf{0.9289}} & \multicolumn{1}{c|}{\textbf{0.9846}} & \multicolumn{1}{c|}{0.9811}   & \multicolumn{1}{c|}{0.9996}   & 1        \\ \hline
\multicolumn{1}{|c|}{$S_0$}                                                                         & \multicolumn{1}{c|}{0.2880}          & \multicolumn{1}{c|}{\textbf{0.5057}} & \multicolumn{1}{c|}{0.6545}          & \multicolumn{1}{c|}{0.9470}   & \multicolumn{1}{c|}{0.9750}   & 0.9756   \\ \hline
\multicolumn{1}{|c|}{$S_1$}                                                                         & \multicolumn{1}{c|}{0.2791}          & \multicolumn{1}{c|}{0.4883}          & \multicolumn{1}{c|}{0.6658}          & \multicolumn{1}{c|}{0.9669}   & \multicolumn{1}{c|}{0.9996}   & 1        \\ \hline
\multicolumn{1}{|c|}{$S_2$}                                                                         & \multicolumn{1}{c|}{0.2790}          & \multicolumn{1}{c|}{0.4968}          & \multicolumn{1}{c|}{\textbf{0.6730}} & \multicolumn{1}{c|}{0.9609}   & \multicolumn{1}{c|}{0.9996}   & 1        \\ \hline
\multicolumn{1}{|c|}{$S_3$}                                                                         & \multicolumn{1}{c|}{0.2897}          & \multicolumn{1}{c|}{0.5030}          & \multicolumn{1}{c|}{0.6574}          & \multicolumn{1}{c|}{0.9151}   & \multicolumn{1}{c|}{0.9716}   & 0.9676   \\ \hline
\multicolumn{1}{|c|}{$S_4$}                                                                         & \multicolumn{1}{c|}{0.1352}          & \multicolumn{1}{c|}{0.1424}          & \multicolumn{1}{c|}{0.1850}          & \multicolumn{1}{c|}{0.2438}   & \multicolumn{1}{c|}{0.1411}   & 0.0647   \\ \hline
\multicolumn{1}{|c|}{$S_5$}                                                                         & \multicolumn{1}{c|}{\textbf{0.2960}} & \multicolumn{1}{c|}{0.3913}          & \multicolumn{1}{c|}{0.5261}          & \multicolumn{1}{c|}{0.9538}   & \multicolumn{1}{c|}{0.9917}   & 0.9968   \\ \hline
\multicolumn{1}{|c|}{$S_6$}                                                                         & \multicolumn{1}{c|}{0.2632}          & \multicolumn{1}{c|}{0.4542}          & \multicolumn{1}{c|}{0.6416}          & \multicolumn{1}{c|}{0.9644}   & \multicolumn{1}{c|}{0.9996}   & 1        \\ \hline
\multicolumn{1}{|c|}{$\Omega$}                                                                      & \multicolumn{1}{c|}{\textbf{0.4433}} & \multicolumn{1}{c|}{\textbf{0.7013}} & \multicolumn{1}{c|}{\textbf{0.8530}} & \multicolumn{1}{c|}{0.9796}   & \multicolumn{1}{c|}{1}        & 1        \\ \hline
\end{tabular}
\caption{LMFDB with top conductor range.}
\label{tab:LMFDB-range}
\end{table}

\begin{figure}
    \centering
    \resizebox{140mm}{!}{\includegraphics{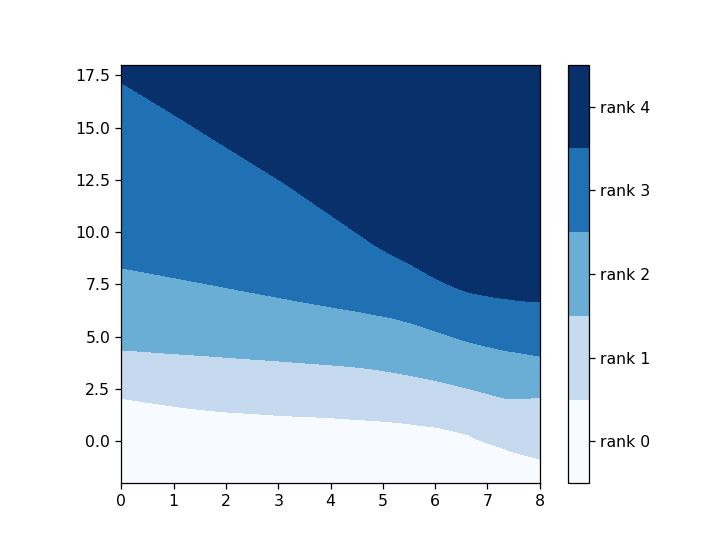}}
    \caption{Rank cutoffs of classifier $S_5$ as a function of a conductor, trained on LMFDB with uniform test set and $p<10^4$. On x-axis are $\log_{10}$ values of conductors and on y-axis are values of sum $S_5$. The unexpected shape of the cutoff between ranks $3$ and $4$ is the consequence of a small number of rank $4$ curves with a small conductor, which are present in the dataset.}
    \label{fig:cutoffs}
\end{figure}

\subsection{Custom dataset}

Test results are presented in Tables \ref{tab:triv_5-uniform} and \ref{tab:triv_5-range}.

In Table \ref{tab:triv_5-uniform} are maximal values of MCC obtained for different classifiers. Classifiers are trained using $80\%$ of randomly selected curves from the custom dataset with a conductor less than $10^{29}$. As in the previous section, MCC is computed on the other $20\%$ of the dataset, not seen during the training.

\begin{table}
\centering
\begin{tabular}{c|cccccc|}
\cline{2-7}
                                                                                                    & \multicolumn{6}{c|}{Number of $a_p$-s used}                                                                                                                                                   \\ \hline
\multicolumn{1}{|c|}{\multirow{2}{*}{\begin{tabular}[c]{@{}c@{}}Type of\\ classifier\end{tabular}}} & \multicolumn{3}{c|}{classify all ranks}                                                                            & \multicolumn{3}{c|}{binary classification}                               \\ \cline{2-7} 
\multicolumn{1}{|c|}{}                                                                              & \multicolumn{1}{c|}{$p<10^3$}        & \multicolumn{1}{c|}{$p<10^4$}        & \multicolumn{1}{c|}{$p<10^5$}        & \multicolumn{1}{c|}{$p<10^3$} & \multicolumn{1}{c|}{$p<10^4$} & $p<10^5$ \\ \hline
\multicolumn{1}{|c|}{CNN}                                                                           & \multicolumn{1}{c|}{\textbf{0.6129}} & \multicolumn{1}{c|}{\textbf{0.7218}} & \multicolumn{1}{c|}{\textbf{0.7958}} & \multicolumn{1}{c|}{0.6695}   & \multicolumn{1}{c|}{0.8335}   & 0.9425   \\ \hline
\multicolumn{1}{|c|}{$S_0$}                                                                         & \multicolumn{1}{c|}{0.5738}          & \multicolumn{1}{c|}{0.6782}          & \multicolumn{1}{c|}{0.7462}          & \multicolumn{1}{c|}{0.6761}   & \multicolumn{1}{c|}{0.8275}   & 0.9175   \\ \hline
\multicolumn{1}{|c|}{$S_1$}                                                                         & \multicolumn{1}{c|}{\textbf{0.5780}} & \multicolumn{1}{c|}{\textbf{0.6890}} & \multicolumn{1}{c|}{\textbf{0.7592}} & \multicolumn{1}{c|}{0.6484}   & \multicolumn{1}{c|}{0.8317}   & 0.9309   \\ \hline
\multicolumn{1}{|c|}{$S_2$}                                                                         & \multicolumn{1}{c|}{0.5649}          & \multicolumn{1}{c|}{0.6761}          & \multicolumn{1}{c|}{0.7521}          & \multicolumn{1}{c|}{0.6407}   & \multicolumn{1}{c|}{0.8252}   & 0.9206   \\ \hline
\multicolumn{1}{|c|}{$S_3$}                                                                         & \multicolumn{1}{c|}{0.5551}          & \multicolumn{1}{c|}{0.6616}          & \multicolumn{1}{c|}{0.7361}          & \multicolumn{1}{c|}{0.6515}   & \multicolumn{1}{c|}{0.8118}   & 0.9074   \\ \hline
\multicolumn{1}{|c|}{$S_4$}                                                                         & \multicolumn{1}{c|}{0.2893}          & \multicolumn{1}{c|}{0.2472}          & \multicolumn{1}{c|}{0.2251}          & \multicolumn{1}{c|}{0.4271}   & \multicolumn{1}{c|}{0.4179}   & 0.3874   \\ \hline
\multicolumn{1}{|c|}{$S_5$}                                                                         & \multicolumn{1}{c|}{0.4987}          & \multicolumn{1}{c|}{0.5990}          & \multicolumn{1}{c|}{0.6696}          & \multicolumn{1}{c|}{0.5230}   & \multicolumn{1}{c|}{0.6919}   & 0.7956   \\ \hline
\multicolumn{1}{|c|}{$S_6$}                                                                         & \multicolumn{1}{c|}{0.5230}          & \multicolumn{1}{c|}{0.6509}          & \multicolumn{1}{c|}{0.7361}          & \multicolumn{1}{c|}{0.5179}   & \multicolumn{1}{c|}{0.7554}   & 0.8961   \\ \hline
\multicolumn{1}{|c|}{$\Omega$}                                                                      & \multicolumn{1}{c|}{\textbf{0.5999}} & \multicolumn{1}{c|}{\textbf{0.7069}} & \multicolumn{1}{c|}{\textbf{0.7807}} & \multicolumn{1}{c|}{0.6821}   & \multicolumn{1}{c|}{0.8527}   & 0.9412   \\ \hline
\end{tabular}

\caption{Custom dataset with uniform test set}
\label{tab:triv_5-uniform}
\end{table}

\begin{figure}
    \centering
    \resizebox{160mm}{!}{\includegraphics{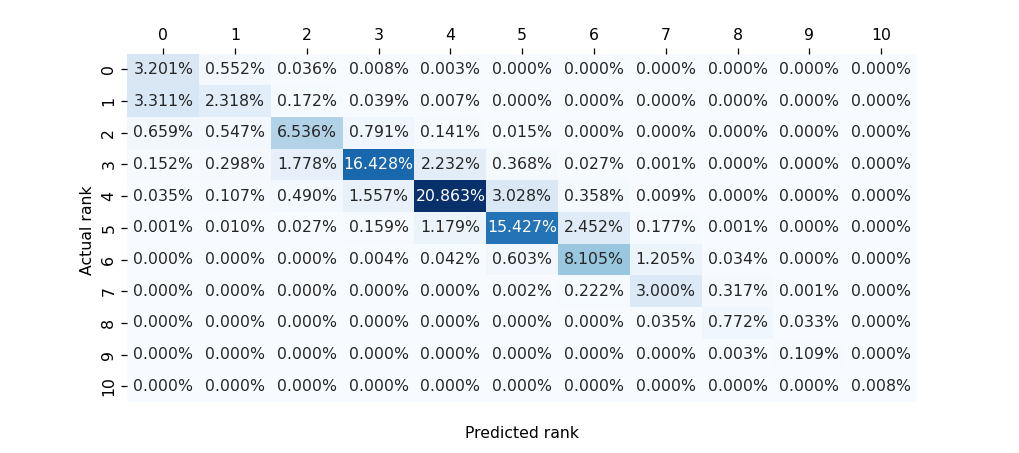}}
    \caption{Confusion matrix of CNN for custom dataset and $p<10^4$ with uniform test set.}
    \label{fig:CM_triv_5_uniform_NN_p10000}
\end{figure}

\begin{figure}
    \centering
    \resizebox{160mm}{!}{\includegraphics{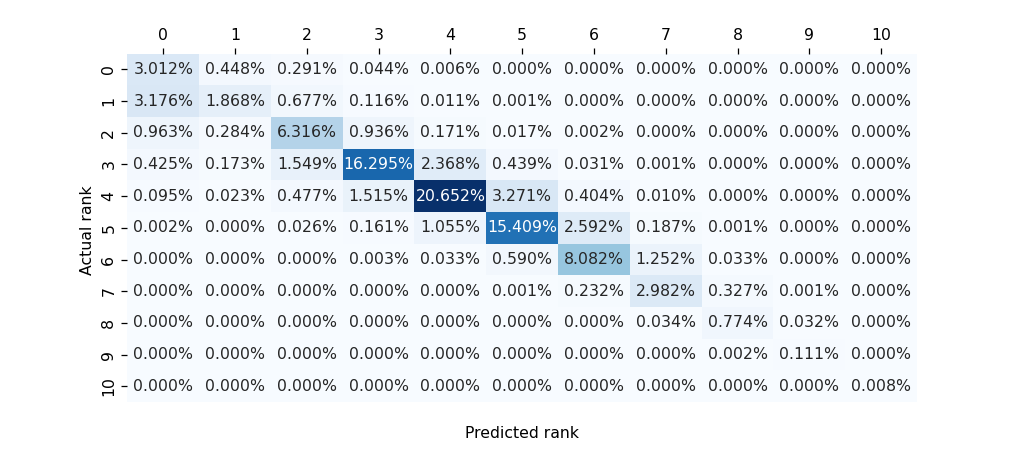}}
    \caption{Confusion matrix of $\Omega$ for custom dataset and $p<10^4$ with uniform test set.}
    \label{fig:CM_triv_5_uniform_ALL_p10000}
\end{figure}

In Table \ref{tab:triv_5-range} classifiers are trained using all curves from the custom dataset with conductor less than $10^{29}$. MCC is computed on all curves from the custom dataset with conductor between $10^{29}$ and $10^{30}$, not seen during the training.

\begin{table}
\centering
\begin{tabular}{c|cccccc|}
\cline{2-7}
                                                                                                    & \multicolumn{6}{c|}{Number of $a_p$-s used}                                                                                                                                                   \\ \hline
\multicolumn{1}{|c|}{\multirow{2}{*}{\begin{tabular}[c]{@{}c@{}}Type of\\ classifier\end{tabular}}} & \multicolumn{3}{c|}{classify all ranks}                                                                            & \multicolumn{3}{c|}{binary classification}                               \\ \cline{2-7} 
\multicolumn{1}{|c|}{}                                                                              & \multicolumn{1}{c|}{$p<10^3$}        & \multicolumn{1}{c|}{$p<10^4$}        & \multicolumn{1}{c|}{$p<10^5$}        & \multicolumn{1}{c|}{$p<10^3$} & \multicolumn{1}{c|}{$p<10^4$} & $p<10^5$ \\ \hline
\multicolumn{1}{|c|}{CNN}                                                                           & \multicolumn{1}{c|}{0.2147}          & \multicolumn{1}{c|}{0.3019}          & \multicolumn{1}{c|}{0.3655}          & \multicolumn{1}{c|}{0.5449}   & \multicolumn{1}{c|}{0.6774}   & 0.8082   \\ \hline
\multicolumn{1}{|c|}{$S_0$}                                                                         & \multicolumn{1}{c|}{0.2533}          & \multicolumn{1}{c|}{\textbf{0.3233}} & \multicolumn{1}{c|}{\textbf{0.3719}} & \multicolumn{1}{c|}{0.5866}   & \multicolumn{1}{c|}{0.6978}   & 0.7900   \\ \hline
\multicolumn{1}{|c|}{$S_1$}                                                                         & \multicolumn{1}{c|}{\textbf{0.2573}} & \multicolumn{1}{c|}{\textbf{0.3291}} & \multicolumn{1}{c|}{\textbf{0.3834}} & \multicolumn{1}{c|}{0.5649}   & \multicolumn{1}{c|}{0.6925}   & 0.7946   \\ \hline
\multicolumn{1}{|c|}{$S_2$}                                                                         & \multicolumn{1}{c|}{0.2340}          & \multicolumn{1}{c|}{0.3118}          & \multicolumn{1}{c|}{0.3688}          & \multicolumn{1}{c|}{0.5472}   & \multicolumn{1}{c|}{0.6814}   & 0.7856   \\ \hline
\multicolumn{1}{|c|}{$S_3$}                                                                         & \multicolumn{1}{c|}{\textbf{0.2556}} & \multicolumn{1}{c|}{0.3189}          & \multicolumn{1}{c|}{0.3645}          & \multicolumn{1}{c|}{0.5688}   & \multicolumn{1}{c|}{0.6830}   & 0.7807   \\ \hline
\multicolumn{1}{|c|}{$S_4$}                                                                         & \multicolumn{1}{c|}{0.1234}          & \multicolumn{1}{c|}{0.1228}          & \multicolumn{1}{c|}{0.1024}          & \multicolumn{1}{c|}{0.3966}   & \multicolumn{1}{c|}{0.3935}   & 0.3842   \\ \hline
\multicolumn{1}{|c|}{$S_5$}                                                                         & \multicolumn{1}{c|}{0.2081}          & \multicolumn{1}{c|}{0.2858}          & \multicolumn{1}{c|}{0.3380}          & \multicolumn{1}{c|}{0.4860}   & \multicolumn{1}{c|}{0.6071}   & 0.6880   \\ \hline
\multicolumn{1}{|c|}{$S_6$}                                                                         & \multicolumn{1}{c|}{0.1803}          & \multicolumn{1}{c|}{0.2757}          & \multicolumn{1}{c|}{0.3527}          & \multicolumn{1}{c|}{0.4321}   & \multicolumn{1}{c|}{0.6217}   & 0.7326   \\ \hline
\multicolumn{1}{|c|}{$\Omega$}                                                                      & \multicolumn{1}{c|}{\textbf{0.2622}} & \multicolumn{1}{c|}{\textbf{0.3246}} & \multicolumn{1}{c|}{\textbf{0.3905}} & \multicolumn{1}{c|}{0.5931}   & \multicolumn{1}{c|}{0.7091}   & 0.8118   \\ \hline
\end{tabular}
\caption{Custom dataset with top conductor range.}
\label{tab:triv_5-range}
\end{table}

\subsection{Discussion}
As expected, classifiers based on $S_0$, $S_1$, and $S_2$ have very similar performance across all tested regimes (improved convergence of sums $S_1$ and $S_2$  does not seem to improve the quality of the classification). They have the best performance of all considered Mestre-Nagao sums. Classifiers based on $S_3$ have overall slightly worse performance than the first three sums. Peculiarly, $S_4$ underperforms in all tests, and its performance, in contrast to every other sum, decreases as the number of used $a_p$-s increases. Classifier $S_5$ performs well but compared to the first four sums is in most cases significantly lesser. While classifier $S_6$ is among the best on the curves having small conductors, its performance fades on larger conductors. Classifier $\Omega$, which is based on all Mestre-Nagao sums simultaneously, is practically always better than any other sum-based classifier, especially on curves of a small conductor.

Finally, CNN based classifiers on the LMFDB dataset were significantly better than $\Omega$ (and consequently any other sum-based classifier). On custom datasets with a uniform test set, they were slightly better than $\Omega$ and other classifiers, but on the custom dataset with top range conductor, they performed worse. Notice that every CNN-based classifier was selected as the best one out of only a $10$ trained model, while each Mestre-Nagao sum-based classifier was selected out of at least $200$ trained models, whose hyperparameters were additionally optimized using the Bayesian optimization technique. We expect that the quality of CNN-based classification would increase with more trained models and possibly surpass other models. As CNN model training is much more computationally expensive, more training was unfeasible for us. It is fascinating how CNN almost perfectly classifies LMDFB with uniform test database (see Figure \ref{fig:CM_LMFDB_uniform_NN_p100000}) - it misclassified only $0.045\%$ of all curves! On the other hand for custom dataset both CNN and $\Omega$ classifiers perform similarly - their confusion matrices look almost identical (see Figures \ref{fig:CM_triv_5_uniform_NN_p10000} and \ref{fig:CM_triv_5_uniform_ALL_p10000}).

Additionally, we noticed that all classifiers made the most mistakes in classifying curves of rank $0$ and $1$. As a byproduct of training Mestre-Nagao sum models, we learned the optimal cutoffs for rank classification (see Figure \ref{fig:cutoffs}), which is of interest when searching for curves of high rank (see Section 7 of \cite{Elkies_new_rank_records}).

\section{Future research}

This paper left us with some open questions which may be addressed in the future projects:

\begin{enumerate}
\item CNN and $\Omega$ models could be used as a substitution for the current use of the Mestre-Nagao sums in the search for elliptic curves of high rank. Also, these models could be trained on other elliptic curves datasets, such as curves with non-trivial torsion, and with additional input of the root number (to employ the Parity conjecture).
\item It is unclear why CNN works much better than all of the other Mestre-Nagao sum based models on the curves from the LMFDB (or in general on curves of a small conductor). Did CNN discover some new mathematics?
\item Suboptimal choices in CNN model architecture and optimizer hyperparameters (see Section \ref{sec:arhitektura_mreze}) may lead to underperforming classifiers. One could perform more extensive optimization of hyperparameters, especially in the case of the custom dataset.
\item Using some of the more advanced neural network architectures, originally developed for natural language processing tasks (see \cite{NIPS2017_Transformer}) may lead to better performing models. The idea is that such architectures can more easily extract distant correlations in the $a_p$-s sequence.
\end{enumerate}

\section*{Acknowledgments}
The first author was supported by the Croatian Science Foundation under the project no.~IP-2018-01-1313, and by the QuantiXLie Center of Excellence, a project co-financed by the Croatian Government and European Union through the European Regional Development Fund - the Competitiveness and Cohesion Operational Programme (Grant KK.01.1.1.01.0004). The second author was supported by Croatian Science Foundation (HRZZ) grant PZS-2019-02-3055 from ``Research Cooperability'' program funded by the European Social Fund.

\bibliographystyle{alpha}
\bibliography{bibliography}
\end{document}